\pdfoutput=1
\documentclass[11pt,letterpaper]{article}

\usepackage[margin=1in]{geometry}
\usepackage{amsmath,amssymb,amsthm,mathtools}
\usepackage[dvipsnames]{xcolor}
\usepackage[colorlinks=true,linkcolor=blue!60!black,citecolor=blue!60!black,urlcolor=blue!60!black]{hyperref}
\usepackage[round]{natbib}
\usepackage[capitalize,nameinlink]{cleveref}
\usepackage{booktabs}

\theoremstyle{plain}
\newtheorem{theorem}{Theorem}
\newtheorem{proposition}[theorem]{Proposition}

\newtheorem{lemma}{Lemma}
\theoremstyle{definition}
\newtheorem{definition}{Definition}
\theoremstyle{remark}
\newtheorem{remark}{Remark}

\crefname{lemma}{Lemma}{Lemmas}
\crefname{theorem}{Theorem}{Theorems}
\crefname{proposition}{Proposition}{Propositions}
\crefname{remark}{Remark}{Remarks}

\newcommand{\cX}{\mathcal{X}}

\newcommand{\cB}{\mathcal{B}}
\newcommand{\E}{\mathbb{E}}
\newcommand{\N}{\mathbb{N}}
\newcommand{\1}{\mathbf{1}}
\newcommand{\phat}{\hat p}
\newcommand{\Bin}{\operatorname{Bin}}
\renewcommand{\Pr}{\operatorname{Pr}}
\DeclarePairedDelimiter{\abs}{\lvert}{\rvert}

\title{Good--Turing for neighborhoods: sharp leave-one-out risk bounds\\
for the mass near values seen $k$ times}
\hypersetup{pdftitle={Good-Turing for neighborhoods: sharp leave-one-out risk bounds for the mass near values seen k times},pdfauthor={Jinran Yi}}
\author{Jinran Yi}
\date{}

\begin{document}
\maketitle

\begin{abstract}
The Good--Turing estimator answers a basic question about an i.i.d.\ sample of size $n$: how much
probability do the values seen exactly $k$ times carry? Often the relevant event concerns similarity
rather than identity: a new observation that is a near-duplicate of, or close to, a value seen exactly
$k$ times. We study the probability of the union of the neighborhoods $N(x)$ of
the values $x$ seen exactly $k$ times, for an arbitrary fixed neighborhood map $N$ on a countable domain.
Its natural estimator is the leave-one-out estimator, which reduces to the Good--Turing estimator for
identity neighborhoods $N(x)=\{x\}$. Because neighborhoods can overlap, the target is in general not a sum
over values. A recent deletion-stability analysis gives mean squared error $O((k+\log n)/n)$ and
leaves open whether the rate $O(k/n)$ holds for every neighborhood map. We show that for all
$1\le k\le n$, every distribution and every neighborhood map, the leave-one-out estimator has mean squared
error at most $53\,g$, where $g=(k/n)/\sqrt{1+k(1-k/n)}$, and that already for identity neighborhoods the
supremum of the risk over distributions is at least $g/9$. Overlapping neighborhoods
therefore cost nothing in order beyond the classical case: the worst-case risk over distributions and
neighborhood maps is of order $\sqrt k/n$ for $1\le k\le n/2$, in particular $O(k/n)$, and
$o(k/n)$ when $k\to\infty$. The logarithm in the stability analysis is nevertheless genuine: for each
fixed $k$, one fixed neighborhood map and a sequence of uniform distributions, the deletion-stability
coefficient is of order $(\log n)/n$ along a sequence of sample sizes, so the rate $O(k/n)$ cannot be
obtained by bounding this coefficient. Finally, the minimax risk for the least favorable domain and
neighborhood map is of order $1/n$ at $k=1$ and of order $1/n^2$ at $k=n$, where the leave-one-out
estimator can have risk one.
\end{abstract}


\section{Introduction}\label{sec:intro}

Given $n$ independent draws from an unknown distribution $D$ on a countable set $\cX$, how much
probability do the values seen exactly $k$ times carry? The Good--Turing estimator answers this question
\citep{good1953population,mcallester2000convergence}. In many settings, however, the event of interest is
similarity rather than identity: whether a new observation is a near-duplicate of, a variant of, or close
to a value seen exactly $k$ times. Such questions concern neighborhoods. Fix a neighborhood map
$N\colon\cX\to2^\cX$, which may encode any notion of similarity. For an integer $k\ge1$, let $c_k(S)$ be
the set of values that occur exactly $k$ times in the sample $S=(x_1,\dots,x_n)$, and let
\[
G_k(S)=\bigcup_{x\in c_k(S)}N(x),\qquad \theta(S)=D(G_k(S)).
\]
Thus $\theta(S)$ is the probability that a fresh draw falls into the neighborhood of some value seen
exactly $k$ times. It is a random quantity: the same sample selects the event and must be used to
estimate its probability. With $N(x)=\{x\}$ it is the total mass of the values seen exactly $k$ times. A
related similarity-based quantity, the mass outside the union of metric balls around all sample points,
is studied by \citet{maurer2022concentration}.

The natural estimator holds out one point at a time: for each $i$ it computes the event $G_k(S^{-i})$ on
the sample $S^{-i}$ with the $i$th point removed and checks whether $x_i$ falls in it,
\[
\phat(S)=\frac1n\sum_{i=1}^n\1\{x_i\in G_k(S^{-i})\}.
\]
For identity neighborhoods $N(x)=\{x\}$ we have $G_k=c_k$, and $x_i\in c_k(S^{-i})$ exactly when $x_i$
occurs $k+1$ times in $S$, so $\phat(S)=(k+1)\Phi_{k+1}(S)/n$, where $\Phi_j(S)$ is the number of values
occurring exactly $j$ times; this is the Good--Turing estimator. In that case the target is a sum over values, and for $0<2k<n$
\citet[Corollary~7]{painsky2022convergence} bounds the mean squared error by a quantity of order
$\sqrt k/n$. With overlapping neighborhoods the target is in general no longer a sum over values. A region can be
covered by several values seen $k$ times or by exactly one, and deleting or adding a single observation
can move a large region into or out of the event.

\citet{feldman2026unified} introduced a general framework for estimating the probability of
data-dependent events and analyzed this leave-one-out estimator through deletion stability. Their
Theorem~5.1 bounds its mean squared error by $1/(2n)$ plus three times their deletion-stability
coefficient, the expected mass of $G_k(S)\,\Delta\,G_k(S^{-n})$, and their Theorem~5.7 bounds this
coefficient by $(3k+2\ln n+2)/n$ for every neighborhood map and every $n\ge2$, which gives the rate
$O((k+\log n)/n)$. They left open whether the rate $O(k/n)$ holds for every neighborhood map
\citep[Section~5.1]{feldman2026unified}. \Cref{thm:risk} answers
this question affirmatively, and \cref{thm:sharp} determines the worst-case order up to absolute
constants.

\paragraph{Results.}
\Cref{thm:risk} shows that for every countable domain, every distribution, every neighborhood map, every
$n\ge2$ and $1\le k\le n$, the mean squared error of $\phat$ is at most
$B_{n,k}=n^{-1}\bigl(\sqrt{6k+9/4}+\sqrt{2k+2+k^2/n}\bigr)^2\le 53k/(2n)$. No symmetry, reflexivity,
transitivity, bounded-overlap or degree assumption on $N$ is needed. \Cref{thm:sharp} sharpens the
dependence on $k$: over all
distributions and neighborhood maps the worst-case risk lies between $\frac19g(n,k)$ and $53\,g(n,k)$,
where $g(n,k)=(k/n)/\sqrt{1+k(1-k/n)}$, and the lower bound already holds for the supremum over distributions with identity neighborhoods
$N(x)=\{x\}$.
Overlapping neighborhoods therefore cost nothing in order beyond the identity case: for $1\le k\le n/2$ the
worst-case risk is $\Theta(\sqrt k/n)$, which is $o(k/n)$ when $k\to\infty$. \Cref{thm:stability} shows
that the logarithm in the stability analysis is nevertheless real: for the single neighborhood map
$N(x)=\cX$ and uniform distributions along a coupon-collector sequence of sample sizes, the
deletion-stability coefficient is asymptotic to $(\log n)/(en)$ for each fixed $k$. Hence the bound of
their Theorem~5.7 is of the right order in $n$ for fixed $k$ along these sample sizes, and no improvement
of it can yield \cref{thm:risk} through their Theorem~5.1. In the same example the target can be read from the sample
without knowing $D$, so the stability coefficient is not a lower bound on the difficulty of estimation.
\Cref{thm:minimax1} shows that the worst-map minimax risk at $k=1$ (the minimax risk for the least
favorable domain and neighborhood map, defined in \cref{sec:setting}) is of order $1/n$, so the rate of
\cref{thm:risk} cannot be improved in $n$ at $k=1$. At the endpoint $k=n$, \cref{prop:endpoint} shows that
the leave-one-out estimator can have risk one, while the worst-map minimax risk is between $1/(10n^2)$
and $4/(n+1)^2$. The proof of \cref{thm:sharp} also bounds the risk by the expectation of a statistic computed from
$\Phi_k$, $\Phi_{k+1}$ and $\Phi_{k+2}$ alone (\cref{rem:observable}).

\paragraph{The idea.}
Write the error as a centered part and a shift. In the centered part, each leave-one-out term
$\1\{x_i\in G_k(S^{-i})\}$ is centered at its own deleted target $D(G_k(S^{-i}))$, so the terms have
conditional mean zero. To bound the covariance of two such terms we condition on the sample with both
points $i$ and $j$ removed. Adding one value can only create the event at some points $z$ (if its
count rises to $k$) or only destroy it (if its count rises above $k$), never both. The regions destroyed
by different added values are disjoint, because each is the part of a neighborhood that no other
value occurring exactly $k$ times covers. Simultaneous exits therefore cost a sum of squared atom masses, whose
expectation is $O(k/n)$ by an occupancy identity, instead of a first-order stability coefficient, which
carries the logarithm. The shift between the average deleted target and the full target has an exact
formula, and its second moment is bounded by a two-deletion switching argument.

\paragraph{What is not claimed.}
We do not claim sharp constants. Apart from $k=1$ (\cref{thm:minimax1}) and $k=n$
(\cref{prop:endpoint}) we do not determine the order of the worst-map minimax risk: for $2\le k\le n-1$
we give only the upper bound of \cref{thm:risk}. \Cref{thm:sharp} concerns only the leave-one-out
estimator, and its lower bound takes a supremum over
distributions with arbitrarily large support.
All our bounds concern the mean squared error for a single sample that both selects the event and is
used to estimate it; beyond Markov's inequality (\cref{sec:related}) we give no high-probability bounds.
The neighborhood map is fixed in advance.
In the neighborhood-correcting learner discussed in Section~8 of \citet{feldman2026unified},
the map $N$ is also fixed, but the learner selects an additional error region using the same sample.
Our results concern $G_k(S)$ alone. \Cref{thm:stability} shows only that \cref{thm:risk}
cannot be proved by bounding the first-order deletion-stability coefficient by $O(k/n)$, since no such
bound holds; it does not exclude other, for example higher-order, stability arguments.

\Cref{sec:setting} fixes notation, recalls the known bounds and states the results.
\Cref{sec:risk,sec:stability,sec:minimax,sec:sharp} contain the proofs, and \cref{sec:related} discusses
related work and open problems.

\section{Setting and main results}\label{sec:setting}

\subsection{Setting}

Throughout, $\cX$ is a countable set, finite or infinite, $D$ is a probability distribution on $\cX$, and
$p_x=D(\{x\})$. For a bounded real-valued function $f$ on $\cX$ we write $Df=\sum_{x\in\cX}p_xf(x)$. A sample
$S=(x_1,\dots,x_n)$ consists of $n$ independent draws from $D$, and following
\citet{feldman2026unified} we write $S^{-i}$ for the sample $x_1,\dots,x_{i-1},x_{i+1},\dots,x_n$ with
the $i$th point removed. All logarithms are natural; we write $\ln$ where we quote bounds of
\citet{feldman2026unified}. For a finite sample or multiset $T$ we write $C_x(T)$ for the multiplicity of
$x$ in $T$, $T+y$ for $T$ with one more copy of $y$, and
\[
c_k(T)=\{x\in\cX: C_x(T)=k\}.
\]
This is the set $c_k(S)$ of \citet[Section~5.1]{feldman2026unified}; we use only $k\ge1$.

Let $N\colon\cX\to2^\cX$ be any fixed function. We make no symmetry, reflexivity, transitivity, bounded-overlap
or degree assumption. Write $N^{-1}(z)=\{x\in\cX: z\in N(x)\}$ for the set of values whose neighborhood contains
$z$. For a multiset $T$ let
\[
G_k(T)=\bigcup_{x\in c_k(T)}N(x),\qquad f_T=\1_{G_k(T)} .
\]
Since $G_k$ depends on $T$ only through its multiplicities, it is a symmetric mapping in the sense of
\citet{feldman2026unified}. For $J\subseteq\cX$ and $j\ge0$ let $F_j(J,T)=\abs{c_j(T)\cap J}$, and let
$\Phi_j(T)=F_j(\cX,T)$ be the number of values that occur exactly $j$ times in $T$. In this notation
$f_T(z)=\1\{F_k(N^{-1}(z),T)>0\}$. On an infinite alphabet $F_0(J,T)$ may be infinite; only the counts
$F_j$ and $\Phi_j$ with $j\ge1$, which are finite, are used below.

\begin{definition}[Target, estimator, risk]\label{def:estimator}
The target is $\theta(S)=D(G_k(S))=Df_S$. The leave-one-out estimator is
\[
\phat(S)=\frac1n\sum_{i=1}^n\1\{x_i\in G_k(S^{-i})\}=\frac1n\sum_{i=1}^n f_{S^{-i}}(x_i),
\]
and its risk is $R_{n,k}(D,N)=\E_{S\sim D^n}\bigl(\phat(S)-\theta(S)\bigr)^2$.
\end{definition}

\paragraph{Known bounds.}
The estimator $\phat$ is the leave-one-out estimator of \citet[Section~5]{feldman2026unified} for the
mapping $S\mapsto G_k(S)$, and $\theta(S)$ is their target. Their deletion-stability coefficient specializes to
\begin{equation}\label{eq:gamma}
\gamma^D(n)=\E_{S\sim D^n}D\bigl(\abs{f_S-f_{S^{-n}}}\bigr)
=\E_{S\sim D^n}D\bigl(G_k(S)\,\Delta\,G_k(S^{-n})\bigr),
\end{equation}
the expected mass that changes when the last point is dropped, and $\gamma(n)=\sup_D\gamma^D(n)$. In this
notation, their Theorem~5.1, stated there for every symmetric mapping into $[0,1]^\cX$, gives
$R_{n,k}(D,N)\le\frac1{2n}+3\gamma^D(n)$ for every $n\ge1$ and every $D$. Their Lemma~5.6 gives
$\gamma(n)\le(2k+1)/n$ for identity neighborhoods $N(x)=\{x\}$, and their Theorem~5.7 gives
$\gamma(n)\le(3k+2\ln n+2)/n$ for every neighborhood map $N$ and every $n\ge2$.
Combining Theorems~5.1 and~5.7 gives, for every $n\ge2$, the explicit bound
\begin{equation}\label{eq:published}
R_{n,k}(D,N)\le\frac1{2n}+3\,\frac{3k+2\ln n+2}{n}=\frac{9k+6\ln n+6.5}{n}.
\end{equation}
This combination is our arithmetic; it is not printed in \citet{feldman2026unified}.

For the minimax statements we use the worst-map minimax risk
\[
R^*_{n,k}=\sup_{\cX,N}\ \inf_{\tilde p}\ \sup_{D}\ \E_{S\sim D^n}\bigl(\tilde p(S)-\theta(S)\bigr)^2,
\]
where the supremum is over countable sets $\cX$ and functions $N\colon\cX\to2^\cX$, and the infimum is
over all estimators $\tilde p$, that is, functions of the sample that may depend on $\cX$, $N$, $n$ and
$k$ but not on $D$.

\subsection{Main results}

\begin{theorem}[Log-free risk]\label{thm:risk}
For every countable $\cX$, every distribution $D$ on $\cX$, every $N\colon\cX\to2^\cX$, every $n\ge2$ and
every integer $1\le k\le n$,
\[
R_{n,k}(D,N)\le B_{n,k}:=\frac1n\Bigl(\sqrt{6k+\tfrac94}+\sqrt{2k+2+\tfrac{k^2}{n}}\Bigr)^2\le\frac{53k}{2n}.
\]
For $k>n$ both the estimator and the target are zero.
\end{theorem}

Hence the leave-one-out estimator attains the rate $O(k/n)$ for every neighborhood map, and the
answer to the question of \citet[Section~5.1]{feldman2026unified} is yes. \Cref{rem:borel} discusses
domains that are not countable. \Cref{tab:compare}
compares $B_{n,k}$ with~\eqref{eq:published}; both are upper bounds on the worst-case risk, not risks.

\begin{table}[ht]
\centering
\begin{tabular}{rrrrr}
\toprule
$n$ & $k$ & $\min\{1,(9k+6\ln n+6.5)/n\}$ & $B_{n,k}$ & ratio\\
\midrule
$10^{3}$ & 1 & 0.05695 & 0.02374 & 2.40 \\
$10^{3}$ & 2 & 0.06595 & 0.03875 & 1.70 \\
$10^{3}$ & 5 & 0.09295 & 0.08366 & 1.11 \\
$10^{4}$ & 1 & 0.00708 & 0.00237 & 2.98 \\
$10^{4}$ & 2 & 0.00798 & 0.00387 & 2.06 \\
$10^{4}$ & 5 & 0.01068 & 0.00836 & 1.28 \\
$10^{5}$ & 1 & 0.00085 & 0.00024 & 3.56 \\
$10^{5}$ & 2 & 0.00094 & 0.00039 & 2.42 \\
$10^{5}$ & 5 & 0.00121 & 0.00084 & 1.44 \\
\bottomrule
\end{tabular}
\caption{The bound~\eqref{eq:published}, obtained by combining Theorems~5.1 and~5.7 of
\citet{feldman2026unified} (our arithmetic), capped at one, compared with $B_{n,k}$ of
\cref{thm:risk}; the last column is their ratio. The bounds are rounded to five decimals; the ratio is
computed from the unrounded values.}
\label{tab:compare}
\end{table}

\begin{theorem}[Stability needs the logarithm]\label{thm:stability}
Let $\cX=\N=\{1,2,\dots\}$ and $N(x)=\cX$ for every $x$. Fix $k\ge1$, let $D_m$ be the uniform
distribution on $\{1,\dots,m\}$, and let
\[
n_m=\bigl\lfloor m\bigl(\log m+k\log\log m-\log k!\bigr)\bigr\rfloor .
\]
Let $F_k(t)$ denote the number of atoms that occur exactly $k$ times among $t$ independent draws from
$D_m$. Then, as $m\to\infty$,
\[
\gamma^{D_m}(n_m)=\frac{k}{n_m}\Pr\{F_k(n_m)=1\}+\frac1m\Pr\{F_k(n_m-1)=1\}
\sim\frac1{em}\sim\frac{\log n_m}{e\,n_m}.
\]
Consequently $\sup_D\gamma^D(n)$ is not $O(k/n)$, and the bound $(3k+2\ln n+2)/n$ of Theorem~5.7 of
\citet{feldman2026unified} is of the right order in $n$ for fixed $k$ along the displayed sequence of sample
sizes. In the same example the target $\theta(S)=\1\{F_k(n)>0\}$ can be read from the sample without
knowing $D$, and the leave-one-out risk satisfies
$R_{n_m,k}(D_m,N)\sim1/(em^2)$, which is of order $(\log n_m)^2/n_m^2$.
\end{theorem}

Thus \cref{thm:risk} cannot be proved by bounding this first-order stability coefficient by $O(k/n)$,
because no such bound holds. This does not exclude other stability-type arguments.

\begin{theorem}[Minimax rate at $k=1$]\label{thm:minimax1}
For every $n\ge2$,
\[
\frac1{4(\sqrt n+1)^2}\le R^*_{n,1}\le B_{n,1}.
\]
Hence $R^*_{n,1}=\Theta(1/n)$. The lower bound holds for one fixed countable $\cX$ and one fixed $N$.
\end{theorem}

\begin{proposition}[Minimax risk at $k=n$]\label{prop:endpoint}
Let $k=n\ge2$.
\begin{enumerate}
\item[(a)] If $D$ is a point mass and $N(x)=\{x\}$ for every $x$, the leave-one-out estimator has risk
$R_{n,n}(D,N)=1$.
\item[(b)] The worst-map minimax risk satisfies $\dfrac{1}{10n^2}\le R^*_{n,n}\le\dfrac{4}{(n+1)^2}$.
\end{enumerate}
\end{proposition}

\begin{theorem}[Sharp frequency dependence]\label{thm:sharp}
Let $g(n,k)=\dfrac{k/n}{\sqrt{1+k(1-k/n)}}$ and let $N_{\mathrm{id}}(x)=\{x\}$. For all integers
$1\le k\le n$,
\[
\frac19\,g(n,k)\le\sup_{\cX,D}R_{n,k}(D,N_{\mathrm{id}})\le\sup_{\cX,D,N}R_{n,k}(D,N)\le\min\{1,\,53\,g(n,k)\},
\]
where the suprema are over countable sets $\cX$, distributions $D$ on $\cX$ and functions
$N\colon\cX\to2^\cX$.
For $1\le k\le n/2$, the worst-case leave-one-out risk is $\Theta(\sqrt k/n)$. This is $o(k/n)$ whenever
$k\to\infty$ in that range. Identity neighborhoods, the classical Good--Turing setting, already give the
displayed lower bound. The constants $1/9$ and $53$ are not optimized.
\end{theorem}

The proof is in \cref{sec:sharp}. The lower bound uses distributions whose support grows without bound;
we do not claim it on a fixed finite alphabet.

\section{Proof of Theorem~\ref{thm:risk}}\label{sec:risk}

If $k>n$ then $c_k(S)=\emptyset$ and $c_k(S^{-i})=\emptyset$ for every $i$, so $\theta(S)=\phat(S)=0$.
Assume from now on that $n\ge2$ and $1\le k\le n$.

\subsection{Two occupancy bounds}

For a multiset $T$ and $s\ge0$ put
\[
Z_s(T)=\sum_{x\in c_s(T)}p_x^2 .
\]

\begin{lemma}[Mass-squared occupancy bound]\label{lem:occupancy}
Let $T$ consist of $m\ge0$ independent draws from $D$. For every integer $k\ge1$,
\begin{equation}\label{eq:occupancy}
\E Z_k(T)\le\frac{k+1}{m+1}.
\end{equation}
\end{lemma}

\begin{proof}
For $k\le m$ and $p\in[0,1]$, since $\binom{m+1}{k+1}=\frac{m+1}{k+1}\binom mk$,
\[
p^2\Pr\{\Bin(m,p)=k\}=\binom mk p^{k+2}(1-p)^{m-k}
=\frac{k+1}{m+1}\,p\,\Pr\{\Bin(m+1,p)=k+1\}.
\]
Hence $\E Z_k(T)=\sum_xp_x^2\Pr\{\Bin(m,p_x)=k\}=\frac{k+1}{m+1}\sum_xp_x\Pr\{\Bin(m+1,p_x)=k+1\}$.
Bounding each probability by one and using $\sum_xp_x=1$ gives~\eqref{eq:occupancy}. If $k>m$ the left
side is zero. All sums have nonnegative terms, so the interchange of sum and expectation is justified for
countable $\cX$.
\end{proof}

For a fixed multiset $T$ and $y\in\cX$, define the entry indicator and the mean entry mass
\[
e_y(z)=\1\{f_T(z)=0,\ f_{T+y}(z)=1\},\qquad a(T)=\E_{Y\sim D}\,De_Y .
\]

\begin{lemma}[Entry probability]\label{lem:entry}
If $T$ consists of $m\ge0$ independent draws from $D$, then
\begin{equation}\label{eq:entry}
\E a(T)\le\frac{k}{m+1}.
\end{equation}
\end{lemma}

\begin{proof}
Fix $z\in\cX$ and let $Y\sim D$ be independent of $T$, so that $T'=T+Y$ is a sample of $m+1$ independent
draws whose last point is $Y$. Deleting a point lowers the multiplicity of its value by one and changes
no other multiplicity. Hence $f_T(z)=0$ and $f_{T'}(z)=1$ hold exactly when $Y$ is the unique element of
$c_k(T')\cap N^{-1}(z)$: any other element of $c_k(T')\cap N^{-1}(z)$ would still lie in
$c_k(T)\cap N^{-1}(z)$, and $Y$ itself has multiplicity $k-1\ne k$ in $T$. If
$c_k(T')\cap N^{-1}(z)$ consists of a single value, exactly $k$ of the $m+1$ positions of $T'$ carry it.
By exchangeability of the $m+1$ positions,
\[
\Pr\{f_T(z)=0,\ f_{T+Y}(z)=1\}=\frac{k}{m+1}\Pr\bigl\{\abs{c_k(T')\cap N^{-1}(z)}=1\bigr\}\le\frac{k}{m+1}.
\]
Integrating over $z\sim D$ gives~\eqref{eq:entry}. No assumption on the overlap of neighborhoods is used.
\end{proof}

\subsection{Exits, entries and a centered covariance bound}

Fix a multiset $T$ and write
\[
d_y(z)=f_{T+y}(z)-f_T(z)=e_y(z)-\ell_y(z),\qquad \ell_y(z)=\1\{f_T(z)=1,\ f_{T+y}(z)=0\}.
\]
Adding $y$ increases only $C_y$. Consequently $e_y$ can be nonzero only when $C_y(T)=k-1$, and $\ell_y$
can be nonzero only when $C_y(T)=k$. For each fixed $y$, $d_y$ is thus an indicator, minus an indicator,
or identically zero. In particular
\begin{equation}\label{eq:unit}
\abs{d_y(z)-Dd_y}\le1,\qquad \abs{\ell_y(z)-D\ell_y}\le1 .
\end{equation}
If $C_y(T)=k$, define the exclusive region
\[
O_y=N(y)\setminus\bigcup_{u\in c_k(T),\,u\ne y}N(u).
\]
Then $\ell_y=\1_{O_y}$, because $c_k(T+y)=c_k(T)\setminus\{y\}$. Call $y$ active if $C_y(T)=k$. The sets
$O_y$ for distinct active $y$ are pairwise disjoint: if $y\ne y'$ are active then $O_{y}$ is disjoint from
$N(y')\supseteq O_{y'}$. This disjointness replaces the logarithmic first-moment stability argument.

For kernels $g_y(x)$ and $h_y(x)$ indexed by $y\in\cX$, define, with $X,Y$ independent draws from $D$,
\begin{equation}\label{eq:form}
\cB(g,h)=\E_{X,Y}\bigl[(g_Y(X)-Dg_Y)(h_X(Y)-Dh_X)\bigr].
\end{equation}
Swapping the roles of $X$ and $Y$ shows that $\cB(g,h)=\cB(h,g)$, and $\cB$ is bilinear.

\begin{lemma}[Conditional centered covariance]\label{lem:cov}
For every multiset $T$,
\begin{equation}\label{eq:cov}
\cB(d,d)\le4a(T)+2Z_k(T).
\end{equation}
\end{lemma}

\begin{proof}
First bound the exit term. Write $b_y=D(O_y)=D\ell_y$, with $b_y=0$ for inactive $y$. Expanding the
product in~\eqref{eq:form} gives
\[
\cB(\ell,\ell)=Q-2H+\beta^2,
\]
where
\[
Q=\sum_{x,y}p_xp_y\,\ell_y(x)\ell_x(y),\qquad
H=\sum_{x,y}p_xp_y\,\ell_y(x)b_x\ge0,\qquad
\beta=\sum_yp_yb_y .
\]
Since the $O_y$ are disjoint, $\sum_y\ell_y(x)\le1$ for each $x$. A nonzero $\ell_x(y)$ also requires
$x\in c_k(T)$. Using $2p_xp_y\le p_x^2+p_y^2$ and the symmetry of the summand in $(x,y)$,
\[
Q\le\sum_xp_x^2\sum_y\ell_y(x)\ell_x(y)\le Z_k(T).
\]
Also $\sum_yb_y\le1$, so the Cauchy--Schwarz inequality gives
\[
\beta^2\le\Bigl(\sum_yb_y\Bigr)\Bigl(\sum_yp_y^2b_y\Bigr)\le Z_k(T),
\]
where the last step uses $b_y\le1$ and $b_y=0$ unless $y\in c_k(T)$. Hence $\cB(\ell,\ell)\le2Z_k(T)$.

Next, since $e_y$ is an indicator,
\[
\E_X\abs{e_y(X)-De_y}=2(De_y)(1-De_y)\le2De_y .
\]
Together with~\eqref{eq:unit}, this implies
\[
\abs{\cB(e,d)}\le2a(T),\qquad\abs{\cB(e,\ell)}\le2a(T).
\]
By bilinearity, symmetry and $d=e-\ell$,
\[
\cB(d,d)=\cB(\ell,\ell)+\cB(e,d)-\cB(e,\ell)\le2Z_k(T)+4a(T).
\]
This proves~\eqref{eq:cov}. The form $\cB$ is not asserted to be positive semidefinite; the proof uses
the explicit expansion and absolute-value bounds only.
\end{proof}

\subsection{The centered leave-one-out residual}

Set
\begin{equation}\label{eq:UV}
A_i=f_{S^{-i}}(x_i)-Df_{S^{-i}},\qquad U=\frac1n\sum_{i=1}^nA_i,\qquad
V=\frac1n\sum_{i=1}^nDf_{S^{-i}}-Df_S .
\end{equation}
Then $\phat(S)-\theta(S)=U+V$ exactly. Conditional on $S^{-i}$, the point $x_i$ is a fresh draw from $D$,
so $f_{S^{-i}}(x_i)$ is a Bernoulli variable with mean $Df_{S^{-i}}$, and
\[
\E A_i^2\le\tfrac14 .
\]
For $i\ne j$, condition on $T=S^{-i,-j}$ (the sample with points $i$ and $j$ removed), let $x=x_i$ and
$y=x_j$, and put $a_0(z)=f_T(z)-Df_T$. Since $S^{-i}=T+y$ and $S^{-j}=T+x$ as multisets,
\[
A_i=a_0(x)+d_y(x)-Dd_y,\qquad A_j=a_0(y)+d_x(y)-Dd_x .
\]
The product $a_0(x)a_0(y)$ has conditional mean zero, since $x$ and $y$ are independent draws from $D$ and
$\E a_0(x)=0$. Each mixed product has conditional mean zero by integrating the variable at which the
corresponding centered function is evaluated: for example, conditional on $x$,
$\E_y[d_x(y)-Dd_x]=0$. Therefore
\[
\E[A_iA_j\mid T]=\cB(d,d),
\]
with $d$ defined relative to $T$. Apply \cref{lem:occupancy,lem:entry,lem:cov} with $m=n-2$:
\[
\E A_iA_j\le\frac{4k+2(k+1)}{n-1}=\frac{6k+2}{n-1}.
\]
There are $n$ diagonal terms and $n(n-1)$ ordered off-diagonal terms in $\E U^2$. Thus
\begin{equation}\label{eq:U}
\E U^2\le\frac1{n^2}\Bigl(\frac n4+n(n-1)\frac{6k+2}{n-1}\Bigr)=\frac{6k+9/4}{n}.
\end{equation}
This controls the centered error without paying the stability coefficient: large exit probabilities
enter only through squared atom masses.

\subsection{Full-sample versus deleted targets}

\begin{lemma}[Deleted-to-full shift]\label{lem:shift}
For every sample $S$ of size $n$ and every $z\in\cX$, with $J=N^{-1}(z)$,
\begin{equation}\label{eq:shift}
\frac1n\sum_{i=1}^nf_{S^{-i}}(z)-f_S(z)
=\frac{k+1}{n}F_{k+1}(J,S)\,\1\{F_k(J,S)=0\}-\frac kn\,\1\{F_k(J,S)=1\}.
\end{equation}
\end{lemma}

\begin{proof}
Recall $f_T(z)=\1\{F_k(J,T)>0\}$, and that deleting a point changes only the multiplicity of its value.
If $F_k(J,S)\ge2$, no deletion removes all elements of $c_k(S)\cap J$, so every term equals $f_S(z)=1$
and both sides vanish. If $F_k(J,S)=1$, deleting one of the $k$ copies of the unique witness lowers
its count to $k-1$ and leaves all other counts unchanged. Any other deletion preserves that witness,
even if it creates a second one. The left side is $(n-k)/n-1=-k/n$. If $F_k(J,S)=0$, then $f_S(z)=0$, and
deleting position $i$ creates an element of $c_k(S^{-i})\cap J$ exactly when $x_i\in J$ and
$C_{x_i}(S)=k+1$; there are $(k+1)F_{k+1}(J,S)$ such positions.
\end{proof}

We need a second-moment estimate for the positive part of~\eqref{eq:shift}, not the first moment that
the stability coefficient controls.

\begin{lemma}[Two-deletion switching bound]\label{lem:twodeletion}
For every fixed $J\subseteq\cX$, every $n\ge2$ and every $k\ge1$,
\begin{equation}\label{eq:twodeletion}
\E\Bigl[\Bigl(\frac{k+1}{n}F_{k+1}(J,S)\,\1\{F_k(J,S)=0\}\Bigr)^2\Bigr]\le\frac{2(k+1)}{n}.
\end{equation}
\end{lemma}

\begin{proof}
Put $r=k+1$ and $W=F_r(J,S)$, and use $W^2=W+W(W-1)$. Each value counted by $W$ occupies $r$ positions,
so $rW\le n$ pointwise, and the diagonal contribution is
$(r/n)^2\E[W\1\{F_k(J,S)=0\}]\le r/n$.

For distinct $u,v\in J$, write $T$ for a sample of $n-2$ independent draws from $D$. We claim
\begin{equation}\label{eq:switch}
\frac{r^2}{n^2}\Pr\{C_u(S)=C_v(S)=r,\ F_k(J,S)=0\}
=\frac{n-1}{n}\,p_up_v\,\Pr\{C_u(T)=C_v(T)=k,\ F_k(J,T)=2\}.
\end{equation}
The condition on the right says that $u$ and $v$ are the only elements of $c_k(T)\cap J$. Indeed, the
probability that a sample of size $\ell$ has the finitely supported multiplicity vector
$(c_x)_{x\in\cX}$ is $\ell!\prod_xp_x^{c_x}/c_x!$. Removing one copy of $u$ and one copy of $v$ is a
bijection from the multiplicity vectors of size $n$ in the event on the left onto those of size $n-2$ in
the event on the right, since the other multiplicities are unchanged. The multinomial formula gives
$r^2$ times the probability of the full vector equal to $n(n-1)p_up_v$ times the probability of the
reduced vector, including when a probability is zero. Summing over the
countably many multiplicity vectors proves~\eqref{eq:switch}.

Sum~\eqref{eq:switch} over ordered pairs of distinct $u,v\in J$; the left side becomes
$(r/n)^2\E[W(W-1)\1\{F_k(J,S)=0\}]$. On a realization of $T$ with exactly two elements $u',v'$ in
$c_k(T)\cap J$, the weighted sum $\sum_{u\ne v}p_up_v\1\{\cdot\}$ equals $2p_{u'}p_{v'}$, which is at
most $p_{u'}^2+p_{v'}^2\le Z_k(T)$; on other realizations it is zero. \Cref{lem:occupancy} with sample
size $n-2$ bounds the off-diagonal contribution by
\[
\frac{n-1}{n}\E Z_k(T)\le\frac{n-1}{n}\cdot\frac{k+1}{n-1}=\frac rn .
\]
Adding the diagonal and off-diagonal contributions proves~\eqref{eq:twodeletion}.
\end{proof}

Integrating~\eqref{eq:shift} against $D(dz)$ gives $V=V_+-V_-$, where
\[
\begin{gathered}
V_+=\sum_zp_z\,\frac{k+1}{n}F_{k+1}(N^{-1}(z),S)\,\1\{F_k(N^{-1}(z),S)=0\},\\
V_-=\sum_zp_z\,\frac kn\,\1\{F_k(N^{-1}(z),S)=1\}.
\end{gathered}
\]
Both are nonnegative, and $V_-\le k/n$. By Jensen's inequality and \cref{lem:twodeletion} applied with
$J=N^{-1}(z)$ for each $z$,
\[
\E V_+^2\le\sum_zp_z\,\E\Bigl[\Bigl(\frac{k+1}{n}F_{k+1}(N^{-1}(z),S)\1\{F_k(N^{-1}(z),S)=0\}\Bigr)^2\Bigr]
\le\frac{2(k+1)}{n}.
\]
Since $(V_+-V_-)^2\le V_+^2+V_-^2$ for nonnegative $V_\pm$,
\begin{equation}\label{eq:V}
\E V^2\le\frac{2(k+1)}{n}+\frac{k^2}{n^2}.
\end{equation}

\subsection{Completion of the proof}

By the triangle inequality in $L^2$, \eqref{eq:U} and~\eqref{eq:V},
\[
\sqrt{R_{n,k}(D,N)}\le\sqrt{\E U^2}+\sqrt{\E V^2}
\le\sqrt{\frac{6k+9/4}{n}}+\sqrt{\frac{2(k+1)}{n}+\frac{k^2}{n^2}},
\]
which is the bound $B_{n,k}$. Using $(s+t)^2\le2s^2+2t^2$ and $k\le n$,
\[
B_{n,k}\le\frac{16k+17/2}{n}+\frac{2k^2}{n^2}\le\frac{18k+17/2}{n}\le\frac{53k}{2n}.
\]
Both $\phat$ and $\theta$ lie in $[0,1]$, so the risk is also at most one. This proves \cref{thm:risk}.
\qed

\begin{remark}[Beyond countable domains]\label{rem:borel}
All results here are stated for countable alphabets. For a possible extension to a standard Borel space
with measurable incidence relation $\{(x,z): z\in N(x)\}$, repeated values arise only at atoms almost
surely. Averages over fresh queries still use the full law $D$. For example, the zero-count probability
for a fresh point would be $\int(1-D(\{x\}))^\ell\,D(dx)$, including the diffuse part. We leave the
general-space formulation outside the scope of this note.
\end{remark}

\section{Proof of Theorem~\ref{thm:stability}}\label{sec:stability}

Let $\cX=\N$, $N(x)=\cX$ for every $x$, and let $D_m$ be uniform on $\{1,\dots,m\}$. For a sample $S$ of
size $t$ from $D_m$, $G_k(S)=\cX$ if $c_k(S)\ne\emptyset$ and $G_k(S)=\emptyset$ otherwise, so
\[
f_S\equiv\1\{F_k(t)>0\},
\]
where $F_k(t)=\Phi_k(S)$ is the number of atoms occurring exactly $k$ times. The sample sizes in
\cref{thm:stability} are in the classical window for collecting $k+1$ copies of every coupon
\citep{erdos1961classical}. Since $D_m(\cX)=1$,
$D_m(G_k(S)\,\Delta\,G_k(S^{-n}))$ is the indicator that $f_S$ and $f_{S^{-n}}$ differ.

\paragraph{The exact formula.}
Let $n\ge2$. Adding the $n$th draw changes only the multiplicity of $x_n$. The event turns on
($F_k(n-1)=0$ and $F_k(n)>0$) exactly when the full sample has exactly one atom of multiplicity $k$ and
$x_n$ is one of its $k$ occurrences. By exchangeability of the $n$ positions, this has probability
$\frac kn\Pr\{F_k(n)=1\}$. The event turns off ($F_k(n-1)>0$ and $F_k(n)=0$) exactly when $S^{-n}$ has
exactly one atom of multiplicity $k$ and $x_n$ equals that atom; since $x_n$ is uniform on
$\{1,\dots,m\}$ and independent of $S^{-n}$, this has probability $\frac1m\Pr\{F_k(n-1)=1\}$. Hence, for
every $m\ge1$ and $n\ge2$,
\begin{equation}\label{eq:gammaexact}
\gamma^{D_m}(n)=\frac kn\Pr\{F_k(n)=1\}+\frac1m\Pr\{F_k(n-1)=1\}.
\end{equation}

\paragraph{Factorial moments.}
The limits $\Pr\{F_k(n_m)=1\}\to e^{-1}$ and $\Pr\{F_k(n_m-1)=1\}\to e^{-1}$ obtained below are instances of
the Poisson limit in \citet[Eqs.~(29)--(30)]{erdos1961classical}. We include the factorial-moment argument
to keep the note self-contained, since the risk calculation also uses variants with fixed numbers of
atoms and draws removed. Fix $k$, let $n_m$ be as in \cref{thm:stability} (for $m$ large enough that $n_m\ge2$), and put
$\lambda_m=n_m/m$. Write $(t)_j=t(t-1)\cdots(t-j+1)$. For each fixed integer $j\ge1$, counting ordered
$j$-tuples of distinct atoms and the placements of exactly $k$ draws on each of them gives
\begin{equation}\label{eq:factorial}
\E(F_k(n_m))_j=(m)_j\,\frac{(n_m)_{kj}}{(k!)^jm^{kj}}\,(1-j/m)^{n_m-kj}\longrightarrow1 .
\end{equation}
The limit follows from $m e^{-\lambda_m}\lambda_m^k/k!\to1$. In detail: with
$\ell_m=\log m+k\log\log m-\log k!$ one has $me^{-\ell_m}(\log m)^k/k!=1$, while
$\ell_m-1/m<\lambda_m\le\ell_m$ and $\lambda_m/\log m\to1$, so $me^{-\lambda_m}\lambda_m^k/k!\to1$. Moreover
$(m)_j\sim m^j$, $(n_m)_{kj}\sim n_m^{kj}$, and
$(1-j/m)^{n_m-kj}=e^{-j\lambda_m}(1+o(1))$ because $n_m/m^2=\lambda_m/m\to0$; hence the right side
of~\eqref{eq:factorial} is $(me^{-\lambda_m}\lambda_m^k/k!)^j(1+o(1))\to1$. The same limit holds with
$n_m$ replaced by $n_m-1$, and with any fixed number of atoms removed and any fixed number of draws
removed, since these changes alter $\lambda_m$ by $O(\lambda_m/m)$ and $m$ by a factor $1+o(1)$.

\paragraph{From moments to probabilities.}
For a nonnegative integer random variable $F$ that is bounded, the finite inclusion--exclusion
expansions are
\[
\Pr(F=0)=\sum_{j\ge0}\frac{(-1)^j}{j!}\E(F)_j,\qquad
\Pr(F=1)=\sum_{j\ge1}\frac{(-1)^{j-1}}{(j-1)!}\E(F)_j .
\]
Their partial sums bracket the probabilities. For the second expansion this follows pointwise from
\[
\sum_{j=1}^L(-1)^{j-1}\frac{(F)_j}{(j-1)!}=F(-1)^{L-1}\binom{F-2}{L-1}\qquad(F\ge2),
\]
while the sum equals $1$ for $F=1$ and $0$ for $F=0$. For the first expansion it follows in the same way
from $\sum_{j=0}^L(-1)^j\binom Fj=(-1)^L\binom{F-1}{L}$ for $F\ge1$, the sum being $1$ for $F=0$. Thus
the partial sums with $L$ odd are upper bounds and those with $L$ even are lower bounds for
$\Pr(F=1)$, and the reverse holds for $\Pr(F=0)$. Let $m\to\infty$ for a fixed truncation $L$
using~\eqref{eq:factorial}, and then let $L\to\infty$. Both probabilities converge to
$\sum_{j\ge0}(-1)^j/j!=e^{-1}$. In particular
\[
\Pr\{F_k(n_m)=1\}\to e^{-1},\qquad\Pr\{F_k(n_m-1)=1\}\to e^{-1}.
\]
By~\eqref{eq:gammaexact}, $\gamma^{D_m}(n_m)\sim1/(em)$, since $k/n_m=o(1/m)$. Finally
$\lambda_m\sim\log m\sim\log n_m$, so $1/m=\lambda_m/n_m\sim\log n_m/n_m$. This proves the stability
statement. Combined with Theorem~5.7 of \citet{feldman2026unified}, the order $(\log n)/n$ is attained
for each fixed $k$ along these sample sizes.

\paragraph{Why this is not a risk lower bound.}
For this $N$, the estimator $\1\{\Phi_k(S)>0\}$ equals the target $\theta(S)$ identically for every
distribution on $\cX$, so its risk is zero. The leave-one-out estimator is different. Since $f_{S^{-i}}$
is constant, $\phat(S)=\frac1n\sum_if_{S^{-i}}(z)$ for any $z$, and \cref{lem:shift} with
$N^{-1}(z)=\cX$ gives, writing $F_j=F_j(n)$,
\[
\phat(S)=\begin{cases}1,&F_k\ge2,\\ 1-k/n,&F_k=1,\\ (k+1)F_{k+1}/n,&F_k=0.\end{cases}
\]
Consequently
\begin{equation}\label{eq:couponrisk}
R_{n,k}(D_m,N)=\frac{k^2}{n^2}\Pr(F_k=1)+\frac{(k+1)^2}{n^2}\E\bigl[F_{k+1}^2\1\{F_k=0\}\bigr].
\end{equation}
Let $r=k+1$ and let $F_k^{(m-j)}(t)$ denote the number of atoms occurring exactly $k$ times among $t$
independent uniform draws from $m-j$ atoms. For $j=1,2$, requiring $j$ specified atoms to have
multiplicity $r$, counting, and then considering the remaining atoms gives
\[
\E\bigl[(F_r(n_m))_j\1\{F_k(n_m)=0\}\bigr]=\frac{(m)_j(n_m)_{jr}}{(r!)^jm^{jr}}(1-j/m)^{n_m-jr}
\Pr\bigl\{F_k^{(m-j)}(n_m-jr)=0\bigr\}.
\]
The prefactor is asymptotic to $(\lambda_m/r)^j$, and the final probability tends to $e^{-1}$ by the
preceding argument. The $j=2$ term dominates the $j=1$ term because $\lambda_m\to\infty$. Substitution
into~\eqref{eq:couponrisk}, using $F^2=(F)_2+F$, yields $R_{n_m,k}(D_m,N)\sim\lambda_m^2/(en_m^2)=1/(em^2)$.
This proves all assertions of \cref{thm:stability}.\qed

\section{Minimax bounds}\label{sec:minimax}

\subsection{Proof of Theorem~\ref{thm:minimax1}}

The upper bound is \cref{thm:risk} with $k=1$, applied to the leave-one-out estimator, which is
available for every $\cX$ and $N$.

For the lower bound use the fixed countable domain $\cX=\{0,1\}\times\N$ and the fixed set
$A=\{1\}\times\N$, and set $N(x)=A$ for every $x$. For an integer $M\ge1$ and $q\in[0,1]$, let $D_{q,M}$
be the distribution of $(\xi,\tau)$, where $\xi$ is Bernoulli$(q)$, $\tau$ is uniform on
$\{1,\dots,M\}$, and they are independent. Let $\Omega_M$ be the event that the $n$ observed tags
$\tau_1,\dots,\tau_n$ are distinct. Its probability
\[
h_{n,M}=\frac{(M)_n}{M^n}\longrightarrow1\qquad(M\to\infty)
\]
does not depend on $q$. On $\Omega_M$ every observation $(\xi_i,\tau_i)$ is distinct, so
$\Phi_1(S)=n\ge1$, $G_1(S)=A$, and the target is $D_{q,M}(A)=q$. The tags, including the event
$\Omega_M$, carry no information about $q$ beyond the Bernoulli observations.

The minimax risk for estimating a Bernoulli parameter from $n$ observations is
\[
r_n=\frac1{4(\sqrt n+1)^2}.
\]
To see this, let $K=\sum_i\xi_i$ and $\alpha=\sqrt n/2$. The posterior mean under the Beta$(\alpha,\alpha)$
prior is
\[
\hat q=\frac{K+\alpha}{n+2\alpha}=\frac{K+\sqrt n/2}{n+\sqrt n}.
\]
Its risk at any $q$ is
\[
\frac{nq(1-q)+n(q-1/2)^2}{(n+\sqrt n)^2}=\frac{n/4}{(n+\sqrt n)^2}=r_n .
\]
It is both a Bayes estimator and an equalizer, so its Bayes risk is $r_n$ and no estimator of $q$ from
$\xi_1,\dots,\xi_n$ has Bayes risk, or maximum risk, smaller than $r_n$.

Now fix an estimator $\tilde p$ for the neighborhood problem. For each fixed $M$, the following lower
bound holds even if the estimator is given $M$. Give $q$ the Beta$(\alpha,\alpha)$ prior.
The maximum risk of $\tilde p$ over $\{D_{q,M}:q\in[0,1]\}$ is at least its prior-averaged risk,
which is at least $\E[\1_{\Omega_M}(\tilde p(S)-q)^2]$ because the loss is nonnegative and $\theta=q$ on
$\Omega_M$. The tags are independent of $(q,\xi_1,\dots,\xi_n)$. Conditional on a tag vector in
$\Omega_M$, the map $(\xi_1,\dots,\xi_n)\mapsto\tilde p(S)$ is an estimator of $q$, whose Bayes risk is
at least $r_n$. Averaging over the tags gives a lower bound $h_{n,M}r_n$ on the maximum risk of
$\tilde p$. Letting $M\to\infty$ shows $\sup_D\E(\tilde p-\theta)^2\ge r_n$ for every $\tilde p$, for
this one fixed $\cX$ and $N$. This proves the lower bound.\qed

\subsection{Proof of Proposition~\ref{prop:endpoint}}

\paragraph{(a)} Let $D=\delta_u$ and $N(x)=\{x\}$. The sample consists of $n$ copies of $u$, so
$c_n(S)=\{u\}$ and $\theta(S)=1$. Each $S^{-i}$ consists of $n-1$ copies of $u$, so
$c_n(S^{-i})=\emptyset$ and $\phat(S)=0$. The risk is one.

\paragraph{(b), upper bound.} At $k=n$, $c_n(S)$ is nonempty exactly when all observations are equal to
some $x$, and then $G_n(S)=N(x)$; otherwise $\theta(S)=0$. Consider the estimator that returns zero unless
all observations equal some $x$, in which case it returns $\1\{x\in N(x)\}$. It is exact unless all
observations are equal. On an all-$x$ sample it estimates $D(N(x))$ within $1-p_x$: if $x\in N(x)$ the
error is $1-D(N(x))\le1-p_x$, and otherwise it is $D(N(x))\le1-p_x$. Hence for $n\ge2$ its risk is at
most
\[
\sum_xp_x^n(1-p_x)^2\le\sup_{p\in[0,1]}p^{n-1}(1-p)^2=\frac{4(n-1)^{n-1}}{(n+1)^{n+1}}\le\frac4{(n+1)^2},
\]
where the supremum is attained at $p=(n-1)/(n+1)$. This holds for every $\cX$, $N$ and $D$.

\paragraph{(b), lower bound.} Let $u\ne w$ be two points of $\cX$, and fix $N(u)=\{w\}$ and
$N(w)=\emptyset$. Compare $\delta_u$ with the law $D'$ given by $D'(u)=1-1/n$ and $D'(w)=1/n$. Under
$\delta_u$ the sample is all-$u$ and the target is $\delta_u(\{w\})=0$. Under $D'$ the all-$u$ sample has
probability $q=(1-1/n)^n$, and on it the target is $D'(\{w\})=1/n$. For $x>1$ the derivative of
$x\log(1-1/x)$ is $1/(x-1)-\log(1+1/(x-1))\ge0$, so $q\ge(1-1/2)^2=1/4$. Let $t$ be the value of a given
estimator on the all-$u$ sample. With equal prior weights on the two laws, its Bayes risk is at least
$\frac12t^2+\frac12q(t-1/n)^2$, whose minimum over $t$ is
\[
\frac{q}{2(1+q)n^2}\ge\frac1{10n^2},
\]
since $q/(1+q)$ is increasing in $q$ and $q\ge1/4$. The maximum risk is at least the Bayes risk. This
proves \cref{prop:endpoint}.\qed

\section{Sharp frequency dependence}\label{sec:sharp}

Throughout this section $n\ge2$, $1\le k\le n$, $r=k+1$, $\Phi_j=\Phi_j(S)$, and $U$, $V$ are as
in~\eqref{eq:UV}, so that $\phat-\theta=U+V$. The bound~\eqref{eq:U} charged the entry terms to their
first moment, which costs $k/n$. Here we charge every term to occupancy counts of the full sample; the
point probabilities of the binomial distribution then supply the factor $[1+k(1-k/n)]^{-1/2}$.

\begin{lemma}[Refined second-moment bounds]\label{lem:refined}
\begin{align*}
\E U^2&\le\frac1{4n}+\frac{\E\bigl[k^2\Phi_k+2k^2\1\{\Phi_k\ge2\}+6kr\Phi_r+r(r+1)\Phi_{r+1}\bigr]}{n^2},\\
\E V^2&\le\frac{\E\bigl[k^2\1\{\Phi_k>0\}+r^2\Phi_r+r(r+1)\Phi_{r+1}\bigr]}{n^2}.
\end{align*}
\end{lemma}

\begin{proof}
\emph{Occupancy switches.} Let $T$ be a sample of $n-2$ independent draws from $D$, and for $\ell\ge0$ and
$j\ge0$ let $\mu_{\ell,j}=\sum_xp_x\Pr\{\Bin(\ell,p_x)=j\}$. For $0\le s\le n-2$ use
$n(n-1)\binom{n-2}{s}=(s+1)(s+2)\binom{n}{s+2}$ and for $0\le j\le n-1$ use
$n\binom{n-1}{j}=(j+1)\binom n{j+1}$. Outside these ranges the corresponding two sides below vanish. Thus
\begin{equation}\label{eq:switches}
\E Z_s(T)=\frac{(s+1)(s+2)}{n(n-1)}\E\Phi_{s+2}(S)\quad(s\ge0),\qquad
\mu_{n-1,j}=\frac{j+1}{n}\E\Phi_{j+1}(S)\quad(j\ge0).
\end{equation}
We also use the following exchangeability device. For an event $\mathcal A(T;x,y)$, the probability of
$\mathcal A(T;X,Y)$ with $T$, $X$, $Y$ independent and $\abs T=n-2$ equals
\[
\frac{1}{n(n-1)}\E\,\#\bigl\{(i,j): i\ne j,\ \mathcal A(S^{-i,-j};x_i,x_j)\bigr\},
\]
where $S$ is a sample of size $n$. Equivalently, take a uniform ordered pair of distinct positions
independently of $S$ and distinguish their values. We call the pairs counted here admissible.

\emph{Diagonal and reduction.} As in \cref{sec:risk}, $\E A_i^2\le1/4$ and, for $i\ne j$,
$\E[A_iA_j\mid T]=\cB(d,d)$ with $T=S^{-i,-j}$. Hence
\[
\E U^2\le\frac1{4n}+\frac{n-1}{n}\E\cB(d,d).
\]

\emph{Centered exits.} We show $\cB(\ell,\ell)\le Z_k(T)$ pointwise. For a query $x$, let $\phi(x)$ be the
unique element of $c_k(T)$ whose neighborhood covers $x$, if there is exactly one, and $\phi(x)=\bot$
otherwise, with $p_\bot=0$. Then $\ell_y(x)=\1\{\phi(x)=y\}$ and $b_y=D\ell_y=D\{x:\phi(x)=y\}$. The
nonzero $b_y$ belong to the finite active set $c_k(T)$ and sum to at most one. Expanding as in the proof
of \cref{lem:cov},
\[
\cB(\ell,\ell)=Q-2\sum_xp_xp_{\phi(x)}b_x+\Bigl(\sum_yp_yb_y\Bigr)^2,\qquad
Q=\sum_{x,y}p_xp_y\ell_y(x)\ell_x(y).
\]
Only loops ($\phi(x)=x$) and reciprocal two-cycles ($\phi(x)=y$, $\phi(y)=x$, $x\ne y$) contribute to
$Q$. Since $\phi$ is a function, distinct cycles use disjoint active vertices, and a two-cycle contributes
$2p_xp_y\le p_x^2+p_y^2$; so $Q$ is at most the sum of $p_x^2$ over active vertices lying on a loop or a two-cycle.
Holding $\phi$ and $p$ fixed, view the right side as a convex function of a free vector $b$.
The actual vector lies in the simplex $\{b\ge0,\ \sum_yb_y\le1\}$ on the active set.
It therefore suffices to bound the function at
$b=0$ and at each unit vector $b=\1_{\{z\}}$ with $z$ active. At $b=0$ the value is $Q\le Z_k(T)$. At
$b=\1_{\{z\}}$ the value is $Q-2p_zp_{\phi(z)}+p_z^2$. If $z$ belongs to neither a loop nor a two-cycle,
it contributes nothing to $Q$, even when it lies on a longer cycle. Then
$Q\le Z_k(T)-p_z^2$, which pays for the added $p_z^2$. A loop at $z$ incurs the subtraction $2p_z^2$. In a
two-cycle $\{z,w\}$ the subtraction $2p_zp_w$ removes that cycle's whole contribution to $Q$, leaving
room for $p_z^2$. Thus $\cB(\ell,\ell)\le Z_k(T)$.

\emph{Entries.} Bilinearity and symmetry give $\cB(d,d)=\cB(e,e)-2\cB(e,\ell)+\cB(\ell,\ell)$, with
\begin{align*}
\cB(e,e)&=\E[e_Y(X)e_X(Y)]-2\E[e_Y(X)\,De_X]+a(T)^2,\\
-2\cB(e,\ell)&=-2\E[e_Y(X)\ell_X(Y)]+2\E[e_Y(X)\,D\ell_X]+2\E[\ell_Y(X)\,De_X]-2a(T)\textstyle\sum_yp_yb_y,
\end{align*}
where the expectations are over independent $X,Y\sim D$ given $T$, and we used the symmetry
$\E[De_Y\,\ell_X(Y)]=\E[\ell_Y(X)\,De_X]$. Discarding the nonpositive terms,
\[
\cB(d,d)\le\cB(\ell,\ell)+\E[e_Y(X)e_X(Y)]+a(T)^2+2\E[e_Y(X)\,D\ell_X]+2\E[\ell_Y(X)\,De_X].
\]
We bound the expectation over $T$ of each of the four entry terms.

(i) Mutual entries. Suppose $X\ne Y$ and $e_Y(X)e_X(Y)=1$. Then $C_X(T)=C_Y(T)=k-1$, and in $S=T+X+Y$ both
$X$ and $Y$ have multiplicity $k$, $c_k(S)=c_k(T)\cup\{X,Y\}$, and no element of $c_k(T)$ covers $X$. So
$Y$ is the only element of $c_k(S)$ other than possibly $X$ itself whose neighborhood contains $X$. Given
$S$, each value $x\in c_k(S)$ therefore has at most one possible partner $y$, and each ordered pair of
distinct values uses $k^2$ ordered pairs of positions; this gives at most $k^2\Phi_k(S)$ admissible
ordered positions. For $X=Y$ the final multiplicity is $r$, and each value of multiplicity $r$ gives
$r(r-1)=kr$ ordered pairs of positions, at most $kr\Phi_r(S)$ in total. By the exchangeability device,
\[
\E\,\E[e_Y(X)e_X(Y)\mid T]\le\frac{\E[k^2\Phi_k+kr\Phi_r]}{n(n-1)}.
\]

(ii) Squared entry mass. Since $a(T)=\sum_zp_z\sum_yp_ye_y(z)$, Jensen's inequality gives
$a(T)^2\le\sum_zp_z\bigl(\sum_yp_ye_y(z)\bigr)^2=\sum_zp_z\sum_{y,y'}p_yp_{y'}e_y(z)e_{y'}(z)$. The terms
with $y=y'$ contribute at most $\sum_yp_y^2De_y\le Z_{k-1}(T)$, since $e_y\ne0$ requires $C_y(T)=k-1$. For
$y\ne y'$ and a query $z$ with $e_y(z)e_{y'}(z)=1$, the two entrants become exactly two elements of
$c_k(T+y+y')\cap N^{-1}(z)$. Treating the query as an independent draw and the two entrants as two
distinguished positions in a sample of size $n$, the exchangeability device gives at most
$2k^2\1\{\Phi_k(S)\ge2\}$ admissible ordered positions. By~\eqref{eq:switches} with $s=k-1$,
\[
\E a(T)^2\le\frac{\E[kr\Phi_r+2k^2\1\{\Phi_k\ge2\}]}{n(n-1)}.
\]

(iii) Mixed terms. For the first, add $Y$ to $T$ to form $R'=T+Y$, a sample of size $n-1$. A nonzero
product $e_Y(X)\,D\ell_X$ requires $C_X(T)=k$ and $C_Y(T)=k-1$, hence $X\ne Y$ and $C_X(R')=k$.
Given the multiset $R'$ and the independent $X$, distinguish a uniform position of $R'$.
The event $e_Y(X)=1$ requires its value $Y$ to be the unique element of $c_k(R')\cap N^{-1}(X)$, so at
most $k$ choices of the distinguished position of $Y$ in $R'$ are possible. Dropping the remaining factor
$D\ell_X\le1$ and averaging the condition $C_X(R')=k$ gives $k\mu_{n-1,k}/(n-1)$. For the other mixed term, write
$De_X=\Pr_Z\{e_X(Z)=1\}$ with a fresh query $Z\sim D$, and add $X$ to $T$ to form $R'=T+X$, leaving $Y,Z$
independent of $R'$. On the joint event $\ell_Y(X)e_X(Z)=1$ the core counts are
$C_Y(T)=k$ and $C_X(T)=k-1$. Thus $X\ne Y$ and $C_Y(R')=k$. Given the multiset $R'$ and $Y,Z$,
at most $k$ distinguished positions of $R'$ could have caused entry at $Z$. Averaging
$\1\{C_Y(R')=k\}$ gives the same bound. By~\eqref{eq:switches},
$k\mu_{n-1,k}/(n-1)=kr\E\Phi_r/(n(n-1))$, so the two mixed terms together, with their factors $2$,
contribute at most $4kr\E\Phi_r/(n(n-1))$.

(iv) Exits. By the centered-exit bound and~\eqref{eq:switches} with $s=k$,
$\E\cB(\ell,\ell)\le\E Z_k(T)=r(r+1)\E\Phi_{r+1}/(n(n-1))$.

Adding (i)--(iv), the coefficient of $kr\E\Phi_r$ is $1+1+4=6$, and
\[
\E\cB(d,d)\le\frac{\E\bigl[k^2\Phi_k+2k^2\1\{\Phi_k\ge2\}+6kr\Phi_r+r(r+1)\Phi_{r+1}\bigr]}{n(n-1)} .
\]
The $n(n-1)$ off-diagonal terms of $\E U^2$ have coefficient $1/n^2$, which yields the bound on $\E U^2$.

\emph{Shift of target.} Write $\Delta_S(z)$ for the right side of~\eqref{eq:shift} with $J=N^{-1}(z)$, so
that $V=\sum_zp_z\Delta_S(z)$, and abbreviate $F_j(z)=F_j(N^{-1}(z),S)$. By Jensen's inequality
$V^2\le\sum_zp_z\Delta_S(z)^2$, and the two terms of $\Delta_S(z)$ have disjoint supports. The negative
square contributes at most $\sum_zp_z(k/n)^2\1\{F_k(z)=1\}\le k^2\1\{\Phi_k>0\}/n^2$. Write
$F_r(z)^2=F_r(z)+(F_r(z))_2$. The diagonal part of the positive square contributes at most
$(r/n)^2\sum_zp_zF_r(z)\le r^2\Phi_r/n^2$. For its off-diagonal part, removing one copy each of two distinct
covering values of multiplicity $r$ leaves exactly two elements of $c_k(T)\cap N^{-1}(z)$; in the reverse
switch~\eqref{eq:switch} the two designated values have weight $2p_up_v\le p_u^2+p_v^2$. As in the proof
of \cref{lem:twodeletion}, for each $z$,
$(r/n)^2\E[(F_r(z))_2\1\{F_k(z)=0\}]\le\frac{n-1}{n}\E Z_k(T)$. Switching $Z_k(T)$ back to $S$
by~\eqref{eq:switches} gives
\[
\frac{n-1}{n}\E Z_k(T)=\frac{n-1}{n}\cdot\frac{r(r+1)\E\Phi_{r+1}}{n(n-1)}=\frac{r(r+1)\E\Phi_{r+1}}{n^2}.
\]
This proves the bound on $\E V^2$.
\end{proof}

\begin{remark}[An observable bound]\label{rem:observable}
Define the sample statistics
\[
\begin{gathered}
\widehat W=\frac1{4n}+\frac{k^2\Phi_k+2k^2\1\{\Phi_k\ge2\}+6kr\Phi_r+r(r+1)\Phi_{r+1}}{n^2},\\
\widehat V=\frac{k^2\1\{\Phi_k>0\}+r^2\Phi_r+r(r+1)\Phi_{r+1}}{n^2}.
\end{gathered}
\]
The right sides of \cref{lem:refined} are $\E\widehat W$ and $\E\widehat V$. Since
$(s+t)^2\le\frac32s^2+3t^2$,
\[
\E(\phat-\theta)^2\le\E\Bigl[\tfrac32\widehat W+3\widehat V\Bigr].
\]
A hash table of raw counts gives $\Phi_k,\Phi_{k+1},\Phi_{k+2}$, hence $\widehat W$ and $\widehat V$, in
expected linear time and space, even when neighborhood queries are expensive. Using
$2\cdot\1\{\Phi_k\ge2\}\le\Phi_k$, $k\1\{\Phi_k>0\}\le k\Phi_k$ and $\sum_{j\ge1}j\Phi_j=n$, one gets
$\widehat W\le(6k+1/4)/n$ and $\widehat V\le(k+1)/n$, so $\frac32\widehat W+3\widehat V\le16k/n$ on every
sample. The inequality holds in expectation only: on a given sample $\frac32\widehat W+3\widehat V$ need
not dominate the realized squared error $(\phat-\theta)^2$, so it is not by itself a confidence radius.
\end{remark}

\begin{lemma}[Largest binomial point probability]\label{lem:mode}
For integers $m\ge1$ and $0\le j\le m$ let $\pi_j(m)=\max_{p\in[0,1]}\Pr\{\Bin(m,p)=j\}$. Then
\[
\frac{1}{3\sqrt{1+j(1-j/m)}}\le\pi_j(m)\le\frac{1}{\sqrt{1+j(1-j/m)}} .
\]
\end{lemma}

\begin{proof}
For $j\in\{0,m\}$ we have $\pi_j(m)=1$ and $j(1-j/m)=0$, so both inequalities hold. Let $0<j<m$, so that
$m\ge2$, and put $\sigma^2=j(1-j/m)=j(m-j)/m$. The function $p\mapsto p^j(1-p)^{m-j}$ is maximized at
$p=j/m$, so $\pi_j(m)=\binom mj\,j^j(m-j)^{m-j}/m^m$. By the bounds
$\sqrt{2\pi t}\,(t/e)^te^{1/(12t+1)}<t!<\sqrt{2\pi t}\,(t/e)^te^{1/(12t)}$ of \citet{robbins1955remark},
valid for every integer $t\ge1$ and applied to $m!$, $j!$ and $(m-j)!$,
\[
\pi_j(m)=\frac{e^{\vartheta}}{\sqrt{2\pi\sigma^2}},\qquad
-\frac1{12j}-\frac1{12(m-j)}<\vartheta<\frac1{12m}-\frac1{12j+1}-\frac1{12(m-j)+1}.
\]
The upper limit for $\vartheta$ is negative because $12j+1<12m$, and the lower limit is at least $-1/6$
because $j\ge1$ and $m-j\ge1$. Moreover $\sigma^2\ge(m-1)/m\ge1/2$, since $j(m-j)\ge m-1$ for
$1\le j\le m-1$. Hence $2\pi\sigma^2\ge1+\sigma^2$, which gives the upper bound
$\pi_j(m)\le(2\pi\sigma^2)^{-1/2}\le(1+\sigma^2)^{-1/2}$. For the lower bound,
$\pi_j(m)\ge e^{-1/6}(2\pi\sigma^2)^{-1/2}\ge e^{-1/6}(2\pi)^{-1/2}(1+\sigma^2)^{-1/2}$, and
$e^{-1/6}/\sqrt{2\pi}>0.337>1/3$.
\end{proof}

\begin{proof}[Proof of \cref{thm:sharp}, upper bound]
The risk is at most one because $\phat$ and $\theta$ lie in $[0,1]$; this also covers $n=1$. Let $n\ge2$
and $m=n-1$, and write $\pi_j=\pi_j(m)$ for $0\le j\le m$ and $\pi_j=0$ for $j>m$.
By~\eqref{eq:switches},
\[
\frac{(j+1)\E\Phi_{j+1}}{n}=\mu_{n-1,j}=\sum_xp_x\Pr\{\Bin(m,p_x)=j\}\le\pi_j .
\]
In \cref{lem:refined} use $2\cdot\1\{\Phi_k\ge2\}\le\Phi_k$ and $\1\{\Phi_k>0\}\le\Phi_k$. Then
\[
\E U^2\le\frac1{4n}+\frac{2k\pi_{k-1}+6k\pi_k+r\pi_{k+1}}{n},\qquad
\E V^2\le\frac{k\pi_{k-1}+r\pi_k+r\pi_{k+1}}{n}.
\]
Put $v=k(1-k/n)$, so that $g(n,k)=(k/n)(1+v)^{-1/2}$. We claim that
$\pi_j\le\sqrt3\,(1+v)^{-1/2}$ for each $j\in\{k-1,k,k+1\}$. For $j>m$ this is trivial. For the legal
indices $j\in\{k-1,k,k+1\}\cap[0,m]$, the quantity $j(1-j/m)$ differs from $v$ by at
most two. Indeed, if $k\le m$, the function $t\mapsto t-t^2/m$ has slope of absolute value at most one on
$[0,m]$, so $\abs{j(1-j/m)-k(1-k/m)}\le1$, and $0\le k(1-k/n)-k(1-k/m)=k^2/(n(n-1))\le1$; if $k=n$, the
only legal index is $j=m$, and both quantities vanish. Hence $1+v\le3\bigl(1+j(1-j/m)\bigr)$, and the
claim follows from \cref{lem:mode}. Since $r\le2k$, we obtain
\[
\E U^2\le\frac1{4n}+10\sqrt3\,g(n,k),\qquad\E V^2\le5\sqrt3\,g(n,k).
\]
Finally, $g(n,k)\ge(k/n)/\sqrt{1+k}\ge1/(\sqrt2\,n)$ for $k\ge1$, so $1/(2n)\le g(n,k)/\sqrt2$. Therefore
\[
R_{n,k}(D,N)\le2\E U^2+2\E V^2\le\Bigl(\frac1{\sqrt2}+30\sqrt3\Bigr)g(n,k)\le53\,g(n,k),
\]
since $30\sqrt3+1/\sqrt2<52.7$. This proves the upper bound.
\end{proof}

\begin{proof}[Proof of \cref{thm:sharp}, lower bound]
Write $P_j=\Pr\{\Bin(n,k/n)=j\}$. Since the maximum of $p\mapsto\Pr\{\Bin(n,p)=k\}$ is attained at
$p=k/n$, we have $P_k=\pi_k(n)$, and \cref{lem:mode} gives
\[
P_k\ge\tfrac13\bigl[1+k(1-k/n)\bigr]^{-1/2}.
\]

\emph{Case $n\ge3k$.} Let $\nu=\lfloor n/k\rfloor$ and use $\nu$ atoms $u_1,\dots,u_\nu$ of mass $p=k/n$
each, with $N(u_i)=\{u_i\}$, and put the leftover mass $1-\nu k/n$ on one further point with empty
neighborhood. Write $C_i=C_{u_i}(S)$. Then $G_k(S)=\{u_i:C_i=k\}$, and $x_j\in G_k(S^{-j})$ exactly when
the observation $x_j$ is a distinguished atom whose full multiplicity is $r$. Hence the error is
$\phat-\theta=n^{-1}\sum_{i=1}^\nu\varepsilon_i$ with
\[
\varepsilon_i=r\1\{C_i=r\}-k\1\{C_i=k\},\qquad\E\varepsilon_i=0,\qquad\E\varepsilon_i^2=k(2k+1)P_k,
\]
where we used the binomial ratio $P_r/P_k=k/r$. For two distinct atoms put
$P_{s,t}=\Pr\{C_i=s,\ C_j=t\}$. Exact multinomial ratios give $P_{r,k}/P_{k,k}=k/r$ and
$P_{r,r}/P_{r,k}=[k-k/(n-2k)]/r$, including the value zero when the last configuration is impossible.
Therefore
\[
\E\varepsilon_i\varepsilon_j=-\frac{k^2P_{k,k}}{n-2k}.
\]
Since $P_{k,k}\le P_k$, the absolute total off-diagonal contribution relative to the diagonal is at most
$(\nu-1)k/[(n-2k)(2k+1)]\le2/(2k+1)$, using $(\nu-1)k\le n-k\le2(n-2k)$. It follows that
\[
R_{n,k}\ge\frac{\nu k(2k-1)P_k}{n^2}.
\]
Here $\nu k/n\ge3/4$, because $\lfloor t\rfloor/t\ge3/4$ for $t\ge3$. Hence
$R_{n,k}\ge\frac34(2k-1)P_k/n\ge\frac34(k/n)P_k\ge\frac14g(n,k)$.

\emph{Case $3k>n$.} Use a single distinguished atom $u$ of mass $p=k/n$ with $N(u)=\cX$, all other
neighborhoods empty. On $\{C_u(S)=k\}$ the target is one; the $k$ deletions of a copy of $u$ give
$G_k(S^{-i})=\emptyset$, and the other deletions give $G_k(S^{-i})=\cX$, so the leave-one-out error
equals $-k/n$. Thus $R_{n,k}\ge p^2P_k\ge(k/(3n))P_k\ge\frac19g(n,k)$. This argument also covers $n=1$.

\emph{Identity neighborhoods.} A point mass with identity neighborhoods has risk one at $k=n$
(\cref{prop:endpoint}(a)) and at $k=n-1\ge1$, where the full sample has no value of multiplicity $n-1$
while every $S^{-i}$ does. To obtain both lower bounds with $N(x)=\{x\}$ for every $x$, spread the leftover
mass in the two constructions uniformly over $L$ further atoms with identity neighborhoods, let the
distinguished atoms ($u_1,\dots,u_\nu$, respectively $u$) have identity neighborhoods, and let
$L\to\infty$ at fixed $n$. The event $H$ that no two draws hit the same further atom satisfies
$\Pr(H^c)\le\binom n2/L$ by a union bound over pairs of draws. On $H$ the further atoms have multiplicity at most one, so they
contribute nothing to the leave-one-out estimate (which, for identity neighborhoods, counts values of
multiplicity $r\ge2$), and they contribute at most $n/L$ to the target, even for $k=1$. The distinguished
atoms contribute $r\Phi^{\circ}_r/n-p\Phi^{\circ}_k$, where $\Phi^\circ_j$ counts distinguished atoms of
multiplicity $j$. Denote this contribution by $\varepsilon$. In the first construction its second moment
equals the risk computed in Case $n\ge3k$. In the second, $\varepsilon=-p$ on $\{C_u=k\}$,
so $\E\varepsilon^2\ge p^2P_k$. In both cases $\abs{\varepsilon}\le1$, and on $H$ the full error
is $\varepsilon-\tau$ with $0\le\tau\le n/L$. Therefore
\[
 R_{n,k}(D,N_{\mathrm{id}})\ge\E[\varepsilon^2\1_H]-2n/L
 \ge\E\varepsilon^2-\binom n2/L-2n/L .
\]
The two error terms vanish as $L\to\infty$ at fixed $n$. Taking the supremum gives
$\sup_{\cX,D}R_{n,k}(D,N_{\mathrm{id}})\ge\frac19g(n,k)$.
\end{proof}

\section{Related work and discussion}\label{sec:related}

For identity neighborhoods $N(x)=\{x\}$, the target is the total mass of the values seen exactly $k$
times, and the leave-one-out estimator is the Good--Turing estimator $(k+1)\Phi_{k+1}/n$
\citep{good1953population}, as observed in Remark~5.1 of \citet{feldman2026unified}.
\citet[Section~2 and footnote~1]{mcallester2000convergence} analyze the version with denominator
$n-k$ for $k<n$ and explicitly distinguish the normalization by $n$ used here.
\citet{mcallester2003concentration} proved concentration inequalities for the missing mass.
For identity neighborhoods and $0<2k<n$, \citet[Corollary~7]{painsky2022convergence} gives the
$O(\sqrt k/n)$ upper bound on the mean squared error of the Good--Turing estimator.
\citet[Theorem~30]{drukh2005concentration} obtain the same scale for the variance of the random
target on a uniform vocabulary in their regime. This gives a lower bound when the estimator uses a
separate independent sample. In \cref{thm:sharp} the estimator observes the sample that selects the
target. The lower construction evaluates this estimator's own error, and the upper bound holds for
arbitrary, possibly overlapping neighborhoods.
For fixed $k$ and an unbounded alphabet, \citet[Theorem~1, Eq.~(19)]{painsky2022convergence} also expresses
the identity-neighborhood risk at sample size $n$, without Poissonization, through $\E\Phi_{k+1}^2$,
$\E\Phi_{k+1}$ and $\E\Phi_{k+2}$, up to an $o(1/n)$ remainder.
Bounds on the estimation error in terms of occupancy counts are also classical in the identity setting:
in the Poissonized infinite urn scheme, \citet[Eq.~(5.3) and Proposition~5.5]{benhamou2017concentration} express the variance of the
Good--Turing estimation error for the missing mass through the expected numbers of values seen once and
twice, and derive confidence bounds from observed occupancy counts. \Cref{rem:observable} is a
fixed-$n$ bound of this kind, in expectation, for exact-frequency unions of arbitrary neighborhoods.
\citet{maurer2022concentration} studies the mass outside the union of metric balls around all sample
points. His Theorem~2.1(iii) bounds the variance of the extended Good--Turing estimator minus the
average deleted-sample target by $3/n$. Maurer credits this bound and its two-deletion proof to
Sourav Chatterjee (private communication to Maurer). An exact-frequency union can both gain and lose
coverage when an observation is added.

The generic deleted-estimate bound appears in \citet[Lemma~1]{devroye1979potential}, who prove, for
symmetric classification rules, a bound of the form $1/(2n)$ plus three times a deletion-instability term.
\citet[Lemma~9, Eq.~(9)]{bousquet2002stability} give its bounded-loss formulation.
\citet[Definition~1]{kale2011cross} analyze cross-validation using a mean-square change in loss
when one observation is replaced by an independent draw. Equation~\eqref{eq:gamma} instead compares
sample sizes $n$ and $n-1$. For indicator-valued events squaring this \emph{deletion} difference
leaves its absolute value unchanged.
The present proof uses the generic decomposition into a centered
part and a shift; what is specific is the bound on simultaneous exits by disjoint regions and the
two-deletion bound on the shift.

\citet{feldman2026unified} formulated the data-dependent event problem, applied the classical
stability guarantee to its leave-one-out estimator, and posed the exact-frequency neighborhood question.
\citet{feldman2026amplification} resolve, up to logarithmic factors, a different question of that paper,
the domain size of hard instances: for every mapping on a known domain of size $d$, their randomized
estimator has mean squared error tending to zero whenever $d=o(n^2/\log^2n)$. It combines sample
amplification with distribution learning, and its guarantee depends on the domain size.

\paragraph{Open problems.}
\Cref{thm:minimax1} and \cref{prop:endpoint} determine the order of the worst-map minimax risk
$R^*_{n,k}$ only at $k=1$ and $k=n$; for $2\le k\le n-1$ the note gives only the upper bound of
\cref{thm:risk}, and the minimax rate there, in particular for $k$ growing with $n$, is open. Together
they show that the leave-one-out estimator is not minimax-rate-optimal uniformly in $k$. For $0<\delta<1$, Markov's inequality and \cref{thm:sharp} give
$\Pr\{\abs{\phat-\theta}>\sqrt{\min\{1,53g(n,k)\}/\delta}\}\le\delta$. An open question is whether the
dependence on $\delta$ can be substantially improved at the same risk scale, uniformly in $D$ and $N$, using
a single sample. Frequency windows and thresholds, that is,
neighborhoods of the values whose multiplicity lies in a given set, are not treated here. For the
neighborhood-correcting learner in Section~8 of \citet{feldman2026unified}, the open question concerns
control of the hybrid estimator with a sample-selected error region and a fixed neighborhood map.
Finally, \cref{thm:risk,thm:stability} together show that first-order deletion stability does not characterize the
rate at which the leave-one-out estimator succeeds; a characterization of estimatability, or of
leave-one-out estimatability, that goes beyond first-order stability is open.

\bibliography{refs}

\end{document}